\documentclass[a4paper,10pt]{article}
\usepackage{amsmath,amssymb}
\usepackage[latin1]{inputenc}
\usepackage[english,frenchb]{babel}
\usepackage[dvips]{epsfig}
\usepackage{mathrsfs}

\newcommand{\pd}{\operatorname{Poids}}
\newcommand{\hyp}{\mathcal{HA}}
\newcommand{\hypc}{\mathcal{HAC}}
\newcommand{\PP}{\operatorname{P\!P}}

\newcommand{\fon}{\mathscr{F}}
\newcommand{\zero}{\hat{0}}
\newcommand{\rg}{\operatorname{Rg}}
\newcommand{\ch}{\operatorname{ch}}
\newcommand{\Stab}{\operatorname{Stab}}
\newcommand{\ind}{\operatorname{Ind}}

\newcommand{\ha}{\mathsf{HA}}
\newcommand{\hac}{\mathsf{HAC}}
\newcommand{\hal}{\mathsf{HAL}}
\newcommand{\Wh}{\operatorname{WH}}
\newcommand{\coho}{\operatorname{H}}

\newcommand{\Y}{\mathsf{Y}}
\newcommand{\YC}{\mathsf{YC}}
\newcommand{\sym}{\mathfrak{S}}
\newcommand{\Sym}{\operatorname{Sym}}

\newcommand{\assoc}{\operatorname{Assoc}}
\newcommand{\comm}{\operatorname{Comm}}
\newcommand{\lie}{\operatorname{Lie}}
\newcommand{\pl}{\operatorname{PreLie}}
\newcommand{\perm}{\operatorname{Perm}}
\newcommand{\cycle}{\operatorname{Cyc}}

\newtheorem{theorem}{Théorème}[section] 
\newtheorem{proposition}[theorem]{Proposition} 
\newtheorem{conjecture}[theorem]{Conjecture} 
\newtheorem{corollary}[theorem]{Corollaire} 
\newtheorem{lemma}[theorem]{Lemme}

\newenvironment{proof}
{\textbf{Preuve.}}
{\hfill\rule{2mm}{2mm}}

\title{Hyperarbres, arbres enracinés \\ et partitions pointées} 
\author{F. Chapoton}
\date{\today}

\begin{document}

\maketitle

\begin{abstract}
  On calcule les polynômes caractéristiques des posets des
  hyperarbres. On montre que la série génératrice de ces polynômes
  fait intervenir les hyperarbres cycliques. On donne aussi une
  conjecture pour l'action du groupe symétrique sur l'homologie de
  Whitney de ces posets. Par ailleurs, on montre que le poset des
  partitions pointées de Vallette est équivalent homotopiquement au
  poset des forêts d'arbres enracinés de Pitman. Le thème commun
  implicite à tous ces objets est la combinatoire de l'opérade PreLie.
\end{abstract}

\selectlanguage{english}

\begin{abstract}
  We compute the characteristic polynomials of the posets of
  hypertrees. We show that the generating series of the polynomials
  can be expressed using cyclic hypertrees. We also propose a
  conjecture on the action of the symmetric groups on the homology of
  these posets. On the other hand, we show that Vallette's poset of
  pointed partitions is homotopy equivalent to Pitman's poset of
  forests. The implicit common thema of the article is the
  combinatorics of the PreLie operad.
\end{abstract}

\selectlanguage{frenchb}

\setcounter{section}{-1}

\section{Introduction}

Les hyperarbres sont des objets combinatoires relativement nouveaux,
notamment par rapport aux arbres. Ils ont été introduits par Berge
dans sa générali\-sation de la théorie des graphes aux hypergraphes
\cite{berge}. Plus récemment, ils ont été utilisés pour étudier
certains sous-groupes d'automorphismes du groupe libre
\cite{mccammond}. L'objet de cet article est d'esquisser un rapport
possible avec la théorie des opérades, plus précisément avec l'opérade
anticyclique PreLie.

Cet article comprend deux parties principales. La première étudie le
poset des hyperarbres sur $n$ sommets. Le résultat principal est le
calcul du polynôme caractéristique, suivant une méthode inspirée par
le calcul par McCammond et Meier du nombre de Möbius du poset obtenu
par l'ajout d'un maximum. Cette description fait intervenir la notion
nouvelle d'hyperarbre cyclique.  On propose ensuite une conjecture
décrivant l'action naturelle du groupe symétrique sur l'homologie de
Whitney de ce poset. On montre que si cette conjecture est vérifiée,
alors l'homologie est fortement liée à l'opérade anticyclique
décrivant les algèbres pré-Lie \cite{prelie,anticyclic}.

La seconde partie est plutôt consacrée aux arbres enracinés. On montre
que le poset des partitions pointés (introduit par Vallette en théorie
des opérades) est équivalent par homotopie au poset des forêts
(introduit par Pitman en probabilités). En utilisant les méthodes de
Sundaram, on donne une description de l'action du groupe symétrique
sur l'homologie de Whitney de ces posets. Dans ce contexte, l'opérade
PreLie joue un rôle explicite et mieux compris que dans le cadre des
hyperarbres.

\section{Ordre partiel sur les hyperarbres}

Un \textbf{hypergraphe} sur un ensemble fini de sommets $I$ est un
ensemble non vide de parties de $I$ de cardinal au moins $2$. Ces
parties sont appelées les \textbf{arêtes} de l'hypergraphe. On peut
définir une notion de chemin entre deux sommets dans un hypergraphe :
formellement c'est une suite alternée d'arêtes et de sommets où chaque
arête contient les deux sommets adjacents. On peut donc parler
d'hypergraphe connexe et de cycle dans un hypergraphe.

Un \textbf{hyperarbre} sur un ensemble fini de sommets $I$ est un
hypergraphe connexe sur $I$ qui ne contient pas de cycle. Ceci
entraîne que deux arêtes distinctes se rencontrent en au plus un
sommet.
\begin{figure}
  \begin{center}
    \scalebox{0.4}{\includegraphics{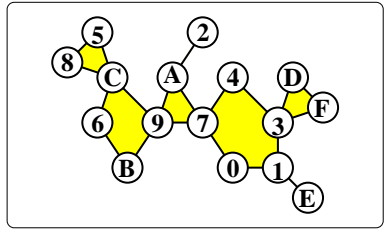}}
    \caption{Un exemple d'hyperarbre sur $\{0,1,2,3,4,5,6,7,8,9,A,B,C,D,E,F\}$.}
    \label{exemple_ha}
  \end{center}
\end{figure}

On va définir un ordre partiel sur l'ensemble des hyperarbres sur $I$.
Ce poset a été étudié par J. McCammond and J. Meier \cite{mccammond}
en relation avec la cohomologie $\ell^2$ de certains groupes
d'automorphismes des groupes libres ; ils montrent notamment que le
poset des hyperarbres est Cohen-Macaulay.

La relation d'ordre sur les hyperarbres est définie comme suit. Un
hyperarbre $S$ est inférieur ou égal à un hyperarbre $T$ si chaque
arête de $S$ est la réunion d'une ou plusieurs arêtes de $T$.

Le rang $\rg(T)$ d'un hyperarbre $T$ est le nombre d'arêtes de $T$
moins un. Le poset des hyperarbres est gradué par le rang. L'unique
élément $\zero$ de rang $0$ est l'hyperarbre dont la seule arête est
l'ensemble $I$ tout entier. Les éléments de rang maximal (égal à
$|I|-2$) sont les $|I|^{|I|-2}$ arbres sur $I$, dont toutes les arêtes
ont cardinal $2$.

\section{Séries génératrices des hyperarbres}

\label{hypera}

On considère ici des séries génératrices classiques pour les
hyperarbres, voir \cite{gessel,mccammond} pour des travaux antérieurs.

Soit $\hyp_n$ l'ensemble des hyperarbres sur l'ensemble
$\{1,\dots,n\}$. On définit le \textbf{poids} d'un hyperarbre $T$
comme le produit
\begin{equation}
  \pd(T)=\frac{1}{t} \prod_{a}(t\, u_{|a|} ),
\end{equation}
où $a$ parcourt les arêtes de $T$. Autrement dit, le poids d'un
hyperarbre est un monôme en les variables $t$ et $(u_i)_{i\geq 2}$. La
puissance de $t$ est le rang de $T$. La puissance de $u_i$ est le
nombre d'arêtes de cardinal $i$ dans $T$.

On introduit alors une série génératrice $\ha$ des hyperarbres selon
le poids, définie comme suit :
\begin{equation}
  \ha=\sum_{n\geq 2} \sum_{T\in\hyp_n} \pd(T) \frac{x^n}{n!}=\frac{x^2}{2} u_2
  +\frac{x^3}{6}(u_3+3 u_2^2)+\dots
\end{equation}

Soit $\ha^p$ la série génératrice similaire des hyperarbres pointés en
un sommet, $\ha^{a}$ celle des hyperarbres pointés en une arête et
$\ha^{pa}$ celle des hyperarbres pointés en un drapeau, c'est-à-dire
munis d'une paire formée d'un sommet distingué et d'une arête
distinguée contenant ce sommet.

On peut facilement obtenir une description récursive de ces séries
génératrices. Tout d'abord, un hyperarbre pointé en un sommet se
décompose naturellement selon les arêtes contenant le sommet
distingué. On obtient ainsi que la série $\ha^{p}$ est caractérisée
par la relation
\begin{equation}
  \label{HAPselonY}
  \ha^p=\frac{x}{t} ( \exp(t\Y) -1 ),
\end{equation}
où la série auxiliaire $\Y$, définie par
\begin{equation}
  \label{YselonHAP}
  \Y=\sum_{n \geq 1} u_{n+1} \frac{(x+t\,\ha^p)^n}{n!},
\end{equation}
est la série génératrice des hyperarbres sur l'ensemble $\{1,\dots,n\}
\sqcup \{\bullet\}$, où le sommet $\bullet$ appartient à une seule
arête.

Par ailleurs, on peut décomposer un hyperarbre pointé en une arête
selon les composantes connexes de l'hypergraphe obtenu en enlevant
cette arête. Chacune de ces composantes est un hyperarbre pointé en
un sommet. On obtient ainsi la relation
\begin{equation}
  \ha^{a}=\sum_{n \geq 1} u_{n+1} \frac{(x+t\,\ha^p)^{n+1}}{(n+1)!}.
\end{equation}

De même, on décrit un hyperarbre pointé en un drapeau en utilisant la
série auxiliaire $\Y$. On obtient l'équation
\begin{equation}
  \ha^{pa}=x \Y \exp(t \Y).
\end{equation}

Enfin, par le principe de dissymétrie (voir \cite[Chap. 4.1]{species}
pour le cas des arbres) qui consiste à utiliser l'existence d'un
centre naturel pour un hyperarbre (qui est soit une arête soit un
sommet), on a
\begin{equation}
  \ha^{pa}+\ha=\ha^{p}+\ha^{a}.
\end{equation}

Toutes ces relations permettent le calcul par récurrence de ces séries
génératrices. On a de plus la relation
\begin{equation}
  \label{derive_ha}
  \ha^p = x\partial_x \ha,
\end{equation}
qui est la traduction habituelle du pointage en un sommet au niveau
des séries génératrices.

\section{Séries génératrices des hyperarbres cycliques}

Un \textbf{hyperarbre cyclique} est un hyperarbre muni en chaque
sommet d'un ordre cyclique sur les arêtes contenant ce sommet.

\label{hyperac}

\subsection{Version complète}

Soit $\hypc_n$ l'ensemble des hyperarbres cycliques sur l'ensemble
$\{1,\dots,n\}$ et soit $\hac$ la série génératrice des hyperarbres
cycliques définie par
\begin{equation}
  \hac=\sum_{n\geq 2} \sum_{T\in\hypc_n} \pd(T) \frac{x^n}{n!}=\frac{x^2}{2} u_2
  +\frac{x^3}{6}(u_3+3 u_2^2)+\dots,
\end{equation}
avec la même définition du poids que précédemment.

Soit $\hac^p$ la série génératrice des hyperarbres cycliques pointés
en un sommet, $\hac^{a}$ celle des hyperarbres cycliques pointés en
une arête et $\hac^{pa}$ celle des hyperarbres cycliques pointés en un
drapeau.

\textbf{Remarque} : on peut aussi voir les hyperarbres cycliques
pointés en un drapeau comme des hyperarbres pointés en un sommet (la
racine) et munis d'un ordre total sur les arêtes entrantes en chaque
sommet (l'arête sortante est l'arête la plus proche de la racine). On
appelle \textbf{hyperarbre ordonné} ce type d'hyperarbre. Pour la
bijection entre ces deux types d'objets, on utilise la structure
arborescente et le drapeau initial pour définir pour chaque sommet
(sauf le sommet pointé) une arête sortante, celle qui est la plus
proche du sommet fixé. Les ordres cycliques sont alors équivalents à
des ordres totaux, en utilisant le drapeau initial (pour la racine) ou
l'arête sortante (pour les autres sommets) pour effectuer la
conversion entre les deux types d'ordres.

\medskip

Comme précédemment pour les hyperarbres, on a une description de ces
séries génératrices par des équations fonctionnelles d'origine
combinatoire.

En utilisant la description par les hyperarbres ordonnés, on voit que
la série $\hac^{pa}$ est caractérisée par la relation
\begin{equation}
  \label{HACPAselonY}
  \hac^{pa}=\frac{x \YC}{1-t \YC},
\end{equation}
où la série auxiliaire $\YC$, définie par
\begin{equation}
  \label{YselonHACPA}
  \YC=\sum_{n \geq 1} u_{n+1} \frac{(x+t\,\hac^{pa})^n}{n!},
\end{equation}
est la série génératrice des hyperarbres ordonnés sur l'ensemble
$\{1,\dots,n\} \sqcup \{\bullet\}$ où le sommet $\bullet$ est la
racine et est contenu dans une seule arête.

Quand on pointe un hyperarbre cyclique en une arête, on obtient un
ensemble (de cardinal au moins $2$) d'hyperarbres ordonnés. On a donc
\begin{equation}
  \hac^{a}=\sum_{n \geq 1} u_{n+1} \frac{(x+t\,\hac^{pa})^{n+1}}{(n+1)!}.
\end{equation}

Quand on pointe en un sommet, on obtient un cycle d'hyperarbres
ordonnés ayant un sommet racine en commun. Ceci entraîne la relation
\begin{equation}
  \hac^{p}=-\frac{x}{t}\ln(1-t Y).
\end{equation}

Enfin, à nouveau par le principe de dissymétrie, on a
\begin{equation}
  \hac^{pa}+\hac=\hac^{p}+\hac^{a}.
\end{equation}

Toutes ces relations permettent le calcul par récurrence de ces séries
génératrices. On a de plus la relation
\begin{equation}
  \hac^p = x\partial_x \hac,
\end{equation}
qui exprime comme précédemment le pointage en un sommet.

\subsection{Version simplifiée}

On spécialise les résultats de la section précédente en posant $u_i=1$
pour tout $i\geq 2$. Le poids d'un hyperarbre devient simplement la variable
$t$ à la puissance la rang de $T$. Par abus de notation, on garde le
même nom pour les séries génératrices, dont les versions complètes ne
seront pas utilisées.

La série $\hac^{pa}$ est déterminée par
\begin{equation}
  \label{def_hac_pa}
  \hac^{pa}=x (\exp(x+t\,\hac^{pa})-1)/(1-t(\exp(x+t\,\hac^{pa})-1)).
\end{equation}

On a la relation suivante :
\begin{equation}
  \label{def_hac_p}
  \hac^p=-\frac{x}{t} \ln(1-t (\exp(x+t\,\hac^{pa})-1)).
\end{equation}

On a aussi
\begin{equation}
  \label{def_hac_a}
  \hac^{a}=\exp(x+t\,\hac^{pa})-1-(x+t\,\hac^{pa}).
\end{equation}

Enfin, on a
\begin{equation}
  \label{def_hac}
  \hac=\hac^{p}+\hac^{a}-\hac^{pa}.
\end{equation}

Bien sûr, on a encore
\begin{equation}
  \label{derive_hac}
  \hac^p = x\partial_x \hac.
\end{equation}

\section{Calcul du polynôme caractéristique}

Dans cette section, on calcule les polynômes caractéristiques des
posets des hyperarbres à l'aide des séries génératrices des
hyperarbres cycliques.

Soit $\chi_n$ le polynôme caractéristique du poset $\hyp_n$, en la
variable $s$ :
\begin{equation}
  \chi_n =\sum_{T \in \hyp_n} \mu(\zero,T) s^{n-2-\rg(T)},
\end{equation}
où $\mu(\zero,T)$ est le nombre de Möbius de l'intervalle $[\zero,T]$.

Soit $\chi_T$ le polynôme caractéristique du poset $\hyp_n^{\geq T}$
formé par les éléments de $\hyp_n$ supérieurs à $T$ :
\begin{equation}
  \chi_T=\sum_{U \geq T} \mu(T,U) s^{n-2-\rg(U)}.
\end{equation}
Comme le poset $\hyp_n^{\geq T}$ est isomorphe au produit $\prod_a
\hyp_{|a|}$ où $a$ parcourt les arêtes de $T$ (voir \cite[Lemme 2.5]{mccammond}),
on a la relation
\begin{equation}
  \chi_T=\prod_a \chi_{|a|}.
\end{equation}
Par conséquent,
\begin{equation}
  \frac{\chi_T}{s} =\frac{1}{s}\prod_a \big{(} s \frac{\chi_{|a|}}{s} \big{)}.
\end{equation}

On sait de plus, par la définition des polynômes caractéristiques, que 
\begin{equation}
  \sum_{T\in\hyp_n} \chi_T=s^{n-2}.
\end{equation}
Donc on a
\begin{equation}
  \sum_{T \in \hyp_n} \frac{1}{s} \chi_T= s^{n-3}.
\end{equation}

On définit des séries $\overline{\ha},\overline{\Y},\dots$ en
remplaçant $u_i$ par $\chi_i /s$ et $t$ par $s$ dans $\ha$, $\Y$, etc.
On obtient alors
\begin{equation}
  \overline{\ha} = \sum_{n \geq 2} s^{n-3} \frac{x^n}{n!}=\frac{1}{s^3}(\exp(s x)-1-s x).
\end{equation}

On en déduit donc, par la relation (\ref{derive_ha}), que
\begin{equation}
  \overline{\ha}^p=\frac{x}{s^2}(\exp(s x)-1).
\end{equation}

Par ailleurs, en inversant la relation (\ref{HAPselonY}), on montre que
\begin{equation}
  s^2 \overline{\Y}=s \ln\left( 1+\frac{s \overline{\ha}^p}{x} \right).
\end{equation}

On a donc obtenu les équations suivantes :
\begin{equation}
  A:=s^2 \overline{\Y}=s \sum_{n\geq 1}\chi_{n+1}
  \frac{(x+s\overline{\ha}^p)^n}{n!}=s \ln\left( 1+\frac{s \overline{\ha}^p}{x} \right),
\end{equation}
et
\begin{equation}
  B:=s x+s^2\,\overline{\ha}^p=s x+ x(\exp(s x)-1).
\end{equation}

Ces relations caractérisent les polynômes $\chi_n$ de la variable $s$.

On effectue ensuite le changement de variables $s=-1/t$ et $x=-t z$.
Soit $\tau_{n+1}$ le polynôme $(-t)^{n-1} \chi_{n+1}(-1/t)$ en la
variable $t$. On a alors les relations suivantes :
\begin{equation}
      A=\sum_{n \geq 1} \tau_{n+1} \frac{B^n}{n!},
\end{equation}
et
\begin{equation}
  \label{probleme}
  \begin{cases}
    A&=-\frac{1}{t}\ln(B/z),\\
    B&=z-t z(\exp(z)-1).
  \end{cases}
\end{equation}
Ces relations caractérisent les polynômes $\tau_n$.

\begin{theorem}
  La solution unique du système {\rm (\ref{probleme})} en fonction
  de $B$ est donnée par
  \begin{equation}
    \begin{cases}
      A&=\partial_B \hac(B),\\
      z&=B+t\,\hac^{pa}(B).
    \end{cases}
  \end{equation}
\end{theorem}

\begin{proof}
  L'unicité est claire. Il suffit donc de vérifier la solution
  proposée. Soient $A$ et $B$ comme dans l'énoncé du théorème. En
  utilisant les relations (\ref{def_hac_pa}),(\ref{def_hac_p}) et
  (\ref{derive_hac}) entre les séries génératrices d'hyperarbres
  cycliques, on a
  \begin{align}
    B A&=\hac^p(B)=-\frac{B}{t} \ln(1-t (\exp(z)-1)),\\
    (z-B)/t&=\hac^{pa}(B)=B (\exp(z)-1)/(1-t(\exp(z)-1)).
  \end{align}

  On obtient en simplifiant
  \begin{align}
    A&=-\frac{1}{t}\ln(1-t(\exp(z)-1)),\\
    B&=z-t z(\exp(z)-1).
  \end{align}
  
  Ceci entraîne que
  \begin{equation}
    A=-\frac{1}{t}\ln(B/z).
  \end{equation}  
\end{proof}

\begin{corollary}
  Le polynôme $\tau_n=(-t)^{n-2}\chi_n(-1/t)$ est la série génératrice
  des hyperarbres cycliques sur $n$ sommets selon le nombre d'arêtes.
  En particulier, le polynôme caractéristique $\chi_n$ du poset
  $\hyp_n$ a des coefficients alternés.
\end{corollary}

La seconde partie du corollaire est déjà connue, car elle résulte du
fait que le poset $\hyp_n$ est Cohen-Macaulay.

Pour retrouver les résultats de McCammond et Meier \cite{mccammond}
sur le nombre de Möbius du poset $\widehat{\hyp}_n$ obtenu en
rajoutant artificiellement un élément maximal à $\hyp_n$, on doit
faire $t=-1$ dans le système (\ref{probleme}), qui se simplifie en
\begin{equation}
  \begin{cases}
    A&=\ln(B/z),\\
    B&=z\exp(z).
  \end{cases}
\end{equation}
La solution est bien connue, donnée par $z=A=W(B)$, où $W$ est la
fonction $W$ de Lambert définie par
\begin{equation}
  W(B)=-\sum_{n \geq 1} n^{n-1} \frac{(-B)^n}{n!}.
\end{equation}

On obtient donc
\begin{equation}
  \mu(\widehat{\hyp}_{n+1})=-\chi_{n+1}(1)=(-1)^n n^{n-1},
\end{equation}
comme attendu.

\section{Homologie de Whitney des hyperarbres}

\subsection{Rappels sur les fonctions symétriques}

Comme référence sur les fonctions symétriques, on renvoie le lecteur
au livre classique de Macdonald \cite{macdonald}. 

On note $\ch$ le caractère de Frobenius qui associe à un module sur le
groupe symétrique une fonction symétrique.

Le pléthysme des fonctions symétriques sera noté $\circ$. La
suspension d'une fonction symétrique $f=f(p_1,p_2,p_3,\dots)$ exprimée
en termes de sommes de puissances est la fonction symétrique $\Sigma_t
f$ définie par
\begin{equation}
  \Sigma_t f = -\frac{1}{t} f(-t p_1,-t^2 p_2, -t^3 p_3,\dots).
\end{equation}
On a clairement $\partial_{p_1} \Sigma_t = (-t) \Sigma_t
\partial_{p_1}$. Par convention, $\Sigma$ désigne la suspension en
$t=1$.

Introduisons quelques fonction symétriques associées à des opérades.
On confond, par abus de notation, une opérade avec la fonction
symétrique associée.

Pour l'opérade associative, on a
\begin{equation}
  \label{def_assoc}
  \assoc=p_1/(1-p_1),
\end{equation}
qui correspond à la somme des représentations régulières des groupes
symétriques.

Soit $\comm$ la fonction symétrique suivante :
\begin{equation}
  \label{def_comm}
  \comm=\exp\left(\sum_{k \geq 1} p_k/k\right)-1.
\end{equation}
Cette fonction symétrique correspond à la somme des représentations triviales.

Le fait que les opérades $\lie$ et $\comm$ soient de Koszul et duales
implique la relation
\begin{equation}
  \label{koszul_comm}
  \Sigma\lie \circ \comm=p_1,
\end{equation}
où $\lie$ est la fonction symétrique associée à l'opérade $\lie$.

On a la relation classique
\begin{equation}
  \label{poisson}
  \comm \circ \lie=\assoc.
\end{equation}

Soit enfin $\perm$ la fonction symétrique définie par
\begin{equation}
  \label{def_perm}
  \perm=p_1 (1+\comm).
\end{equation}
Le fait que les opérades $\pl$ et $\perm$ soient de Koszul et duales
implique la relation
\begin{equation}
  \label{koszul_perm}
  \Sigma\pl \circ \perm= \perm \circ \Sigma\pl =p_1,
\end{equation}
où $\pl$ est la fonction symétrique associée aux arbres enracinés
\cite{prelie}. Par l'interprétation usuelle de l'action de $p_1
\partial_{p_1}$ comme le pointage en un sommet, la fonction $p_1
\partial_{p_1} \pl$ correspond donc aux arbres doublement pointés. Un
tel objet se décompose de façon unique (en utilisant l'unique chemin
joignant les deux points marqués) en une liste d'arbres enracinés. On
a donc la relation
\begin{equation}
  \label{vertebres}
  p_1 \partial_{p_1} \pl = \assoc \circ \pl.
\end{equation}

Des relations (\ref{def_perm}) et (\ref{koszul_perm}), on déduit la
relation
\begin{equation}
  \label{def_pl}
  \Sigma\pl \,(\comm \circ \Sigma\pl)=p_1-\Sigma\pl.
\end{equation}

\begin{lemma}
  \label{somme1}
  On a
  \begin{equation}
    p_1 \partial_{p_1}
    \Sigma\pl + \partial_{p_1} (\comm \circ \Sigma\pl)=1.
  \end{equation}
\end{lemma}

\begin{proof}
  Ceci équivaut à
  \begin{equation}
    \Sigma(p_1 \partial_{p_1}\pl) + \partial_{p_1} (\comm \circ \Sigma\pl)=1.
  \end{equation}
  En utilisant la relation (\ref{vertebres}), le membre de gauche devient
  \begin{equation}
    \Sigma \assoc \circ \Sigma \pl + (\partial_{p_1}\Sigma\pl)(1+\comm \circ \Sigma\pl).
  \end{equation}
  En utilisant l'expression (\ref{def_assoc}) de $\assoc$ et les
  relations (\ref{def_pl}) et (\ref{vertebres}), ceci devient
  \begin{equation}
    \frac{\Sigma\pl}{1+\Sigma\pl}-\Sigma\left(\frac{\pl}{p_1(1-\pl)}\right)\left(\frac{p_1}{\Sigma\pl}\right),
  \end{equation}
  ce qui donne bien $1$.
\end{proof}

\subsection{Quelques fonctions symétriques nouvelles}

\label{nouvelles}

Introduisons de nouvelles fonctions symétriques, inspirées de celles
associées aux hyperarbres cycliques, mais distinctes.

Soit $\hal^{pa}$ la fonction symétrique avec un paramètre $t$ définie par
\begin{equation}
  \label{def_hal_pa}
  \hal^{pa}=p_1 [\Sigma_t \assoc] \circ \comm \circ [p_1+ (-t) \hal^{pa}].
\end{equation}

Notons qu'ici et dans toute la suite, le pléthysme agit de façon non
triviale sur la variable $t$ : on a $p_k \circ t=t^k$. La variable $t$
peut donc seulement être spécialisée en $0$ ou $1$.

Soit $\hal^{p}$ la fonction symétrique avec un paramètre $t$ définie par
\begin{equation}
  \label{def_hal_p}
  \hal^{p}=p_1 [\Sigma_t \lie] \circ \comm \circ [p_1+ (-t) \hal^{pa}].
\end{equation}

Soit $\hal^{a}$ la fonction symétrique avec un paramètre $t$ définie par
\begin{equation}
  \label{def_hal_a}
  \hal^{a}=[\comm - p_1] \circ [p_1+ (-t) \hal^{pa}].
\end{equation}

Soit enfin $\hal$ la fonction symétrique avec un paramètre $t$ définie
par
\begin{equation}
  \label{def_hal}
  \hal=\hal^{p}+\hal^{a}-\hal^{pa}.
\end{equation}

Ces relations permettent de calculer ces fonctions symétriques par
récurrence.

Par les formules (\ref{poisson}),(\ref{def_hal_pa}) et
(\ref{def_hal_p}), on a aussi la relation suivante :
\begin{equation}
  \frac{\hal^{pa}}{p_1}=[\Sigma_t \comm] \circ \frac{\hal^{p}}{p_1}.
\end{equation}

\begin{proposition}
  On a la relation
  \begin{equation}
    \hal^{p}=p_1 \partial_{p_1} \hal.
  \end{equation}
\end{proposition}

\begin{proof}
  On pose $C=\comm \circ (p_1 -t \hal^{pa})$. On calcule d'abord, en
  utilisant la relation (\ref{def_hal_pa}), la
  dérivée $\partial_{p_1}\hal^{pa}$ en fonction de $C$. On obtient
  \begin{equation}
    \partial_{p_1}\hal^{pa}=\frac{C}{1+t C} 
    + p_1 (1-t\partial_{p_1}\hal^{pa}) \frac{1+C}{(1+tC)^2}.
  \end{equation}
  En simplifiant les dénominateurs, ceci donne la relation
  \begin{equation}
    \label{deri_hal_pa}
    (1+2 t C+t^2C^2+t p_1+t p_1 C)\partial_{p_1}\hal^{pa}
    =C+t C^2+p_1+p_1 C.
  \end{equation}

  On utilise ensuite la relation (\ref{def_hal}) pour calculer $p_1
  \partial_{p_1}\hal$. On obtient
  \begin{equation}
    \hal^{p}+p_1^2 (1-t\partial_{p_1}\hal^{pa})\frac{1+C}{1+tC}+p_1
    (1-t\partial_{p_1}\hal^{pa}) C -p_1 \partial_{p_1}\hal^{pa}.
  \end{equation}

  Pour démontrer la proposition, il faut donc vérifier la relation
  \begin{equation}
     \partial_{p_1}\hal^{pa}
    =p_1 (1-t\partial_{p_1}\hal^{pa})\frac{1+C}{1+tC}+
    (1-t\partial_{p_1}\hal^{pa}) C.
  \end{equation}
  Mais ceci résulte immédiatement de la relation (\ref{deri_hal_pa})
  obtenue plus haut.
\end{proof}

\subsection{Description de l'homologie de Whitney}

Soit $P$ un poset gradué avec un élément minimum $\zero$ et dont
toutes les chaînes maximales ont la même longueur. Les groupes
d'homologie de Whitney de $P$, dénotés $\Wh_i(P)$, ont été introduits
par Baclawski \cite{bacla} et étudiés par Björner \cite{bjorn}. Si le
poset $P$ est Cohen-Macaulay, alors on a la relation suivante avec le
polynôme caractéristique :
\begin{equation}
  \sum_{i=0}^r s^{r-i} (-1)^i \dim \Wh_i(P) = \chi_P(s), 
\end{equation}
où $r$ est le rang maximal dans $P$.

Les groupes d'homologie de Whitney ont la description suivante :
\begin{equation}
  \Wh_i(P)=\bigoplus_{x \in P} \widetilde{\coho}_{i-2}(\zero,x),
\end{equation}
où les $\widetilde{\coho}$ sont les groupes d'homologie réduite des
intervalles de $P$.

Soit $\Wh$ la fonction symétrique génératrice des caractères de
l'action des groupes symétriques sur l'homologie de Whitney des posets
$\hyp_n$ :
\begin{equation}
  \Wh:=\sum_{n \geq 2} \sum_{i=0}^{n-2} \ch(\Wh_i(\hyp_n)) (-t)^i.
\end{equation}

On propose la conjecture suivante décrivant l'action du groupe
symétrique sur l'homologie de Whitney du poset des hyperarbres.

\begin{conjecture}
  \label{grosse}
  On a la relation suivante
  \begin{equation}
    \Wh=\hal.
  \end{equation}
\end{conjecture}

Cette conjecture est vraie au niveau des dimensions. Pour le vérifier,
il suffit de comparer l'expression des polynômes caractéristiques
obtenue plus haut en fonction de $\hac$ (en $u_i=1$ pour $i\geq 2$) et
la spécialisation en $p_k=0$ pour $k\geq 2$ des fonctions symétriques
$\hal$.

Cette comparaison utilise notamment le fait que la spécialisation de
la fonction symétrique $\lie$ est $-\ln(1-p_1)$.

\subsection{Caractéristiques d'Euler et relation avec $\pl$}

On substitue $t=1$ dans les formules de la section \ref{nouvelles}.
Les fonctions symétriques obtenues sont notées $\overline{\hal}$, etc.

On sait (\cite[Lemme 1.1]{sunda1}) que si la conjecture \ref{grosse}
est vérifiée, alors $\overline{\hal}$ correspond (aux signes près) à
l'action sur l'homologie du poset $\widehat{\hyp}_n$ obtenu en
ajoutant un maximum artificiel au poset $\hyp_n$.

Soit $\overline{\hal}^{pa}$ la fonction symétrique avec un paramètre
$t$ définie par
\begin{equation}
  \overline{\hal}^{pa}=
  p_1 \left[\frac{p_1}{1+p_1}\right] \circ \comm \circ [p_1 - \overline{\hal}^{pa}].
\end{equation}

Ceci se simplifie en la relation
\begin{equation}
  \overline{\hal}^{pa}=
  (p_1-\overline{\hal}^{pa}) \comm \circ (p_1-\overline{\hal}^{pa}).
\end{equation}

Par comparaison avec la relation (\ref{def_pl}), on a donc
\begin{equation}
  \overline{\hal}^{pa}=p_1-\Sigma\pl.
\end{equation}

On obtient ensuite
\begin{equation}
  \overline{\hal}^{p}=p_1 [\Sigma\lie \circ \comm] \circ [p_1 -
  \overline{\hal}^{pa}]=p_1 [p_1 - \overline{\hal}^{pa}]= p_1 \Sigma \pl.
\end{equation}

On a aussi la relation
\begin{equation}
  \label{deri_hal}
  \frac{1}{p_1}\overline{\hal}^{p}=\partial_{p_1} \overline{\hal}= \Sigma \pl.
\end{equation}

On calcule 
\begin{equation}
  \overline{\hal}^{a}=(\comm-p_1)\circ [p_1 -
  \overline{\hal}^{pa}]=\comm \circ \Sigma\pl -\Sigma\pl.
\end{equation}

Enfin, on obtient
\begin{equation}
  \overline{\hal}=-p_1+\Sigma\pl+p_1 \Sigma\pl + \comm \circ \Sigma\pl 
  -\Sigma\pl,
\end{equation}
ce qui donne
\begin{equation}
  \label{hal_pl}
  \overline{\hal}=-p_1+p_1 \Sigma\pl+ \comm \circ \Sigma\pl.
\end{equation}

Mais on sait, par la relation (\ref{def_pl}), que
\begin{equation}
  \comm \circ \Sigma\pl=-1+p_1/\Sigma\pl.
\end{equation}

Donc on a
\begin{equation}
  1 +\overline{\hal}=p_1(-1+ \Sigma\pl +1/\Sigma\pl).
\end{equation}

On rappelle que l'action des groupes symétriques sur l'opérade
anticyclique $\pl$ est donnée par la fonction symétrique $M$
caractérisée par la relation suivante (voir \cite[Eq.
(50)]{anticyclic}) :
\begin{equation}
  M + 1=p_1 (1+ \pl +1/\pl).
\end{equation}

Comme la suspension est anti-multiplicative,
\begin{equation}
  \Sigma M - 1 =-p_ 1 (-1+ \Sigma\pl +1/\Sigma\pl).
\end{equation}

Donc
\begin{proposition}
  On a la relation
  \begin{equation}
    \overline{\hal}=-\Sigma M.
  \end{equation}  
\end{proposition}

Si la conjecture \ref{grosse} est vérifiée, ceci devrait fournir une
relation entre la suspension de l'opérade anticyclique $\pl$ et
l'homologie du poset $\widehat{\hyp}_n$.

% \subsection{Lien avec Brace}

% Remarque : calculons la composante de degré maximal en $t$ de
% $\hal^{pa}$. On pose $v=1/t$.

% On applique la suspension $\Sigma_v$ à l'équation (\ref{def_hal_pa}) :
% \begin{equation}
%   \Sigma_v \hal^{pa} = - v p_1 \Sigma_v(\Sigma_t\assoc) \circ
%   (\Sigma_v \comm) \circ (p_1 + (-t)\Sigma_v \hal^{pa}).
% \end{equation}
% Soit, en posant $K_v=-\frac{1}{v} \Sigma_v \hal^{pa}$,
% \begin{equation}
%   K_v= p_1 \assoc \circ (\Sigma_v \comm) \circ (p_1+K_v).
% \end{equation}
% On peut alors poser $v=0$ :
% \begin{equation}
%   K_0=p_1 \assoc \circ (p_1 + K_0).
% \end{equation}

% On en déduit que $K_0$ est une somme de Catalan copies de la
% représentation régulière. Ceci est, au terme de degré $1$ près, la
% fonction symétrique associée à l'opérade Brace. On rappelle que
% l'opérade Brace est munie d'un morphisme vers l'opérade $\pl$ par
% symétrisation.

\section{Partitions pointées et forêts d'arbres enracinés}

\subsection{Partitions pointées}

Le poset des partitions pointées $\PP_I$ d'un ensemble $I$, introduit
par Vallette dans \cite{vallette} et étudié ensuite dans
\cite{chapvall}, est une variante intéressante du poset classique des
partitions.

Une partition pointée de $I$ est la donnée d'une partition de $I$ et
d'un élément distingué (``pointé'') dans chaque part de cette
partition. Par abus de notation, on note par la même lettre une
partition pointée et la partition sous-jacente. La relation d'ordre
est la suivante : une partition pointée $\nu$ est inférieure à $\pi$
si la partition $\nu$ raffine la partition $\pi$ (chaque part de $\pi$
est l'union de parts de $\nu$) et si les éléments pointés dans $\pi$
sont aussi pointés dans $\nu$. Le poset des partitions pointés a un
élément minimal donné par la partition en singletons, et $|I|$
éléments maximaux correspondant aux différents pointages de la
partition de $I$ en une seul part.

\subsection{Forêts d'arbres enracinés}

Le poset des forêts $F_I$ a été introduit par Pitman dans
\cite{pitman}, pour des motivations en probabilité. Il est aussi utile
en combinatoire, voir \cite[Ex. 5.29]{stanley2}.

Une forêt d'arbres enracinés sur $I$ est un graphe sur l'ensemble de
sommets $I$ dont les composantes connexes sont simplement connexes
(des arbres) et munies chacune d'un élément distingué appelé la
racine. On peut alors orienter les arêtes vers les racines. Un point
de vue équivalent est de considérer un forêt comme un ensemble
d'arêtes orientées ayant les propriétés adéquates (absence de cycles
et de fourches divergentes).

La relation d'ordre est la suivante : une forêt $G$ est inférieure à
une forêt $F$ si on peut obtenir $G$ en enlevant des arêtes orientés à
$F$. Le poset des forêts a un unique élément minimal donné par la
forêt dont les arbres sont des singletons, et $|I|^{|I|-1}$ éléments
maximaux qui correspondent aux arbres enracinés sur $I$.

Les intervalles dans $F_I$ sont des posets booléens. La fonction de
Möbius est donc $-1$ à la puissance le rang. On voit aussi que le
poset $F_I$ provient d'un complexe simplicial sur l'ensemble des
arêtes orientées. Ceci entraîne que sa réalisation géométrique est la
subdivision barycentrique de ce complexe simplicial.

\subsection{Comparaison}

On a un morphisme de poset $\phi$ de $F_I$ dans $\PP_I$ qui associe à
une forêt la partition de $I$ formée par les arbres et le pointage de
chaque part donné par la racine de chaque arbre. Cette application est
équivariante pour l'action du groupe symétrique.

\begin{theorem}
  \label{homotopie}
  L'application $\phi$ induit une équivalence en homotopie.
\end{theorem}

\begin{proof}
  On utilise le critère suivant : il suffit de vérifier que la fibre
  de chaque élément est contractile (Lemme fibre de Quillen, voir
  \cite{wachs}). On choisit de regarder la fibre vers le bas.

  Comme la fibre d'une partition pointée en plusieurs morceaux est
  isomorphe comme poset au produit des fibres des morceaux, il suffit
  de montrer que la fibre d'une partition pointée en un seul bloc est
  contractile. On fixe donc un ensemble fini $I$ et un élément pointé
  $i$ dans $I$. On note $p_i$ cette partition pointée.

  Par la définition de la relation d'ordre de $F_I$ par enlèvement
  d'arêtes, la réalisation géométrique de $F_I$ est naturellement une
  subdivision barycentrique du complexe simplicial $X$ sur l'ensemble
  $I^2$ (une paire (a,b) est vue comme une arête $(a \leftarrow b)$)
  dont les simplexes sont les ensembles d'arêtes de forêts. Ce
  complexe simplicial est pur et les simplexes maximaux sont en
  bijection avec les arbres enracinés sur $I$.

  La fibre par $\phi$ de la partition pointée $p_i$ correspond alors
  au sous-complexe simplicial $X_1$ dont les simplexes sont les
  ensembles d'arêtes contenant au moins une arête de la forme $(i
  \leftarrow j)$. Les simplexes correspondent aux forêts ayant $i$
  parmi leurs racines et les simplexes maximaux correspondent aux
  arbres ayant $i$ pour racine. Ce complexe simplicial est pur.

  On définit, pour $1\leq r \leq |I|-1$, un complexe simplicial $X_r$
  comme le sous complexe simplicial de $X$ dont les simplexes sont les
  ensembles d'arêtes contenant au moins $r$ arêtes de la forme $(i
  \leftarrow j)$. Ce sont des complexes simpliciaux purs. Les
  simplexes maximaux correspondent aux arbres ayant $i$ pour racine et
  tels que la valence de $i$ est au moins $r$.

  On a des inclusions naturelles $X_{|I|-1} \subset \dots \subset X_1$.

  On va montrer que $X_r$ est contractile par récurrence descendante
  sur $r$.
  
  Considérons le cas $r=|I|-1$. Il y a un seul simplexe maximal,
  correspondant à la corolle de racine $i$, formé par les arêtes $(i
  \leftarrow j)$ pour $j\in I \setminus \{i\}$. Ce complexe est
  évidemment contractile.

  Supposons maintenant que $X_r$ est contractile, pour un certain
  $2\leq r \leq |I|-1$. On va montrer que $X_{r-1}$ se rétracte sur
  $X_r$, donc est aussi contractile.

  Soit $C$ un simplexe maximal de $X_{r-1}$ correspondant à un arbre
  dont la racine $i$ est de valence exactement $r-1$. Ceci signifie
  que $C$ contient exactement $r-1$ arêtes de type $(i \leftarrow j)$
  et donc $|I|-1-(r-1)$ autres arêtes.  Considérons les simplexes de
  codimension $1$ au bord de $C$ qui contiennent $r-2$ arêtes de type
  $(i\leftarrow j)$. Chacun de ces simplexes est contenu dans un seul
  simplexe maximal de $X_{r-1}$, qui est $C$. On peut donc écraser $C$
  sur le reste de son bord. Cette partie restante du bord de $C$ est
  formée de simplexes contenant $r-1$ arêtes de type $(i \leftarrow
  j)$. Ce bord est donc contenu dans le sous complexe $X_{r}$.

  Ceci montre que $X_{r-1}$ se rétracte sur $X_{r}$ qui est
  contractile par récurrence, donc $X_{r-1}$ aussi. Ceci termine la
  récurrence.

  On a donc montré que la fibre $X_1$ est contractile. Ceci termine la
  démonstration.  
\end{proof}

Comme corollaire du théorème \ref{homotopie}, on obtient que les
groupes d'homologie de $\PP_I$ et $F_I$ sont isomorphes comme modules
sur les groupes symétriques, de par l'invariance de l'homologie par
équivalence homotopique \cite[Th. 5.2.2]{wachs}. 

\section{Homologie de Whitney}

\subsection{Calcul pour les partitions pointées}

On sait (\cite[Th. 1.3]{chapvall}) que le polynôme caractéristique du poset
$\PP_I$ est
\begin{equation}
  (s-|I|)^{|I|-1}.
\end{equation}
On vérifie aussi (voir \cite{stanley2,kozlov}) que ce polynôme est
également le polynôme caractéristique du poset $F_I$, par la propriété
que tous les intervalles dans $F_I$ sont des posets booléens et en
utilisant la série génératrice connue $(s+|I|)^{|I|-1}$ pour les
forêts.

On obtient ici une description de l'action du groupe symétrique sur
l'homologie de Whitney des partitions pointées. 

On se place pour ce qui suit dans la catégorie des espaces vectoriels
gradués avec la règle des signes de Koszul (\textit{i.e.} la catégorie
des complexes avec différentielle nulle).

\begin{proposition}
  Le caractère de l'homologie de Whitney du poset des partitions
  pointées est donné par
  \begin{equation}
    \comm \circ \Sigma_t \pl.
  \end{equation}
\end{proposition}

\begin{proof}
  Les orbites du groupe symétrique $\sym_n$ sur $\PP_n$ sont indexées
  par les partitions de $n$. Soit $x_\lambda$ un représentant de
  l'orbite $\lambda$. On note $\lambda_i$ le nombre de parts de taille
  $i$ dans une partition $\lambda$.

  Le stabilisateur de $x_\lambda$ est le sous-groupe
  \begin{equation}
    \Stab(x_\lambda)=\prod_{i} \sym_{\lambda_i} [\sym_{i-1}], 
  \end{equation}
  où $\sym_{\lambda_i}[\sym_{i-1}]$ est un produit en couronne. On
  introduit aussi le groupe un peu plus gros défini par
  \begin{equation}
    G_\lambda=\prod_{i} \sym_{\lambda_i} [\sym_{i}]. 
  \end{equation}

  Alors, on a la description suivante, voir \cite[Th.
  5.1]{bjorn},\cite[Th. 1.2]{sunda1}, de l'action de $\sym_n$ sur
  l'homologie de Whitney :
  \begin{equation}
    \label{whpp}
    \Wh(\PP_n)\simeq \bigoplus_{\lambda} \ind_{\Stab(x_\lambda)}^{\sym_n} \widetilde{\coho}(\zero,x_\lambda).
  \end{equation}

  On décompose l'induction en deux étapes :
  \begin{equation}
    \bigoplus_{\lambda} \ind_{G_\lambda}^{\sym_n}
    \ind_{\Stab(x_\lambda)}^{G_\lambda} \widetilde{\coho}(\zero,x_\lambda).
  \end{equation}
  
  Par l'isomorphisme de Künneth, on décompose l'homologie réduite en
  produit selon la taille des parts de $\lambda$ (voir \cite[Prop.
  2.1]{sunda2}) et on décompose l'induction en produit d'inductions :
  \begin{equation}
    \bigoplus_{\lambda} \ind_{G_\lambda}^{\sym_n}
    \bigotimes_i  \ind_{\sym_{\lambda_i} [\sym_{i-1}]}^{\sym_{\lambda_i} [\sym_{i}]} \widetilde{\coho}(\zero,x_{(i^{\lambda_i})}).
  \end{equation}
  
  Par la description de l'homologie pour les puissances d'un poset
  (voir \cite[Prop. 2.3]{sunda2}), on a
  \begin{equation}
    \widetilde{\coho}(\zero,x_{(i^{\lambda_i})})\simeq \Sym^{\lambda_i}
    \widetilde{\coho}(\zero,x_{(i)}),
  \end{equation}
  où $\Sym$ désigne la puissance symétrique.

  On obtient alors (par compatibilité entre produit en couronnes et
  puissance symétrique)
  \begin{equation}
    \bigoplus_{\lambda} \ind_{G_\lambda}^{\sym_n}
    \bigotimes_i  \Sym^{\lambda_i} \left(\ind_{\sym_{i-1}}^{\sym_{i}} \widetilde{\coho}(\zero,x_{(i)})\right).
  \end{equation}
  
  Or, par un cas particulier (terme de degré maximal) de la formule
  (\ref{whpp}) et les résultats de \cite{vallette} sur l'homologie de
  $\PP_n$, on a
  \begin{equation}
    \operatorname{Sgn}(n)\otimes_{\sym_n}\pl(n)=\Wh_{n-1}(\PP_n)=\ind_{\sym_{n-1}}^{\sym_{n}} \widetilde{\coho}(\zero,x_{(n)}),
  \end{equation}
  où ces espaces sont placés en degré $n-1$ et $\operatorname{Sgn}(n)$
  est la représentation signe.

  Donc on a obtenu
  \begin{equation}
    \bigoplus_{\lambda} \ind_{G_\lambda}^{\sym_n}
    \bigotimes_i  \Sym^{\lambda_i}  \left(\operatorname{Sgn}(i) \otimes_{\sym_i}\pl(i) \right).
  \end{equation}
  
  Au niveau des caractères, on reconnaît dans cette expression le
  pléthysme
  \begin{equation}
    \comm \circ \Sigma_t \pl.
  \end{equation}
\end{proof}

On a aussi une description explicite du caractère de ce module. Par
commodité, on donne plutôt le résultat équivalent pour une suspension.

\begin{proposition}
  \label{CE_prelie}
  La fonction caractéristique de $\Sigma_t \comm \circ \pl$ est
  \begin{equation}
    \sum_{\lambda,|\lambda|\geq 1}(\lambda_1-t)^{\lambda_1-1}\prod_{k \geq
      2}\left( (f_k(\lambda)-t^k)^{\lambda_k}- k \lambda_k (f_k(\lambda)-t^k)^{\lambda_k-1}\right) \frac{p_\lambda}{z_\lambda},
  \end{equation}
  où la somme porte sur les partitions non vides $\lambda$,
  $\lambda_k$ est le nombre de parts de taille $k$ dans la partition
  $\lambda$ et $f_k(\lambda)$ est le nombre de points fixes de la
  puissance $k^{eme}$ d'une permutation de type cyclique $\lambda$.
  Les $p_{\lambda}$ sont les fonctions symétriques sommes de
  puissances et les $z_{\lambda}$ des constantes classiquement
  associées aux partitions.
\end{proposition}

\begin{proof}
  On a la relation
  \begin{equation}
   \Sigma\perm\circ \pl=\left( p_1 \exp\left(-\sum_{k\geq 1} p_k/k\right)\right) \circ \pl =p_1
  \end{equation}

  On introduit de nouvelles variables
  \begin{equation}
    y_\ell= p_\ell \circ \pl.
  \end{equation}

  On a alors la relation inverse
  \begin{equation}
    p_\ell=y_\ell \exp\left(-\sum_k y_{k \ell}/k\right).
  \end{equation}

  Soit $\lambda$ une partition. Pour calculer le coefficient de
  $p_\lambda$ dans la fonction symétrique $\Sigma_t \comm \circ \pl$,
  il faut calculer le résidu
  \begin{equation}
    \iiint(\Sigma_t \comm \circ \pl) \prod_{i=1}^{r} \frac{dp_i}{p_i^{\lambda_i+1}}.
  \end{equation}
  On peut supposer sans restriction que toutes les variables $y_j$ et
  $p_j$ pour $j>r$ sont nulles. La suspension de $\comm$ est
  \begin{equation}
    \Sigma_t \comm=\frac{-1}{t}\left(\exp\left(-\sum_{k\geq 1} t^k p_k/k\right)-1\right).
  \end{equation}
  On effectue le changement de variables pour obtenir une intégrale en
  les variables $y$. On a
  \begin{equation}
    \prod_{i=1}^{r} dp_i
    =\exp\left(-\sum_{i} \sum_k y_{ik}/k\right)\prod_{i=1}^{r} (1-y_i) dy_i. 
  \end{equation}

  A un facteur $-t$ et au terme constant près, on a donc à calculer le résidu
  \begin{equation}
    \iiint \exp\left(-\sum_{k\geq 1} t^k y_k/k\right) \exp\left(\sum_{i}\lambda_i \sum_k y_{ik}/k\right)\prod_{i=1}^{r} \frac{(1-y_i)}{y_i^{\lambda_i+1}} dy_i. 
  \end{equation}

  En regroupant les exponentielles et en inversant les sommations, on
  trouve
  \begin{equation}
    \iiint \exp\left(\sum_{k\geq 1} (f_k(\lambda)-t^k) y_k/k\right) \prod_{i=1}^{r} \frac{(1-y_i)}{y_i^{\lambda_i+1}} dy_i.    
  \end{equation}

  Cette intégrale se décompose en un produit de résidus en chaque
  variable $y_i$. On obtient facilement le résultat attendu.
\end{proof}

\subsection{Relation avec l'opérade $\pl$}

A une fonction symétrique $f=\sum_{n\geq 1} f_n$ correspond une suite
de modules $\fon_n$ sur les groupes symétriques et un foncteur $\fon$ qui
associe à un espace vectoriel $V$ l'espace vectoriel
\begin{equation}
  \fon(V)=\oplus_{n\geq 1} \fon_n\otimes_{\sym_n} V^{\otimes n}.
\end{equation}
Les foncteurs ainsi obtenus sont appelés foncteurs analytiques.
Réciproquement, on peut retrouver la suite de modules $(\fon_n)_n$ à
partir du foncteur $\fon$.

\medskip

La fonction symétrique $\comm \circ \Sigma_t \pl$ sert ainsi à décrire
le complexe de Chevalley-Eilenberg calculant la cohomologie des
algèbres pré-Lie libres vues comme algèbres de Lie. En effet, la
fonction symétrique $\pl$ correspond au foncteur qui associe à $V$
l'algèbre pré-Lie libre sur $V$. La suspension et la composition avec
$\comm$ correspondent à prendre l'algèbre extérieure.

Comme on sait que les algèbres pré-Lie libres sont des algèbres de Lie
libres \cite{foissy}, la cohomologie du complexe de
Chevalley-Eilenberg est concentrée en un seul degré et correspond aux
générateurs. La caractéristique d'Euler
\begin{equation}
  \comm \circ \Sigma \pl
\end{equation}
correspond donc aux générateurs des algèbres pré-Lie libres comme
algèbres de Lie.

Bien sûr, comme dans le cas des hyperarbres, cette caractéristique
d'Euler décrit aussi l'homologie du poset obtenu en rajoutant un
élément maximal au poset des partitions pointées.

On obtient ci-dessous une description de cette fonction symétrique à
l'aide des fonctions symétriques $\overline{\hal}$. Si la conjecture
\ref{grosse} est vraie, ceci donne une relation homologique non
triviale entre le poset des hyperarbres et le poset des partitions
pointées.

\begin{proposition}
  On a la relation suivante :
  \begin{equation}
    \comm \circ \Sigma \pl
    = p_1 - (p_1 \partial_{p_1} \overline{\hal}-\overline{\hal} ).
  \end{equation}
\end{proposition}

\begin{proof}
  On va calculer 
  \begin{equation}
    \comm \circ \Sigma \pl +p_1 \partial_{p_1}
  \overline{\hal}-\overline{\hal}.
  \end{equation}
  En utilisant (\ref{hal_pl}) et (\ref{deri_hal}), on obtient
  \begin{equation}
    \comm \circ \Sigma \pl + p_1 \Sigma\pl +p_1 - p_1
    \Sigma\pl -\comm \circ \Sigma\pl.
  \end{equation}
\end{proof}

\textbf{Remarque :} Les composantes homogènes de $\comm \circ \Sigma_t \pl$
apparaissent aussi implicitement dans les travaux de McCammond et
Meier comme décrivant l'action du groupe symétrique sur l'algèbre de
cohomologie du groupe des automorphismes symétriques du groupe libre
\cite{mac_coho}.

\section{Annexe}

Pour éventuelle référence ultérieure, on rassemble ici les formules
essentielles qui décrivent l'action des groupes symétriques sur les
hyperarbres et les hyperarbres cycliques. Les preuves, essentiellement
basées sur les descriptions combinatoires des sections \ref{hypera} et
\ref{hyperac}, sont omises.

\subsection{Caractère des hyperarbres}

Soit $\ha^{p}$ la fonction symétrique avec un paramètre $t$ définie par
\begin{equation}
  \ha^{p}=p_1 [\Sigma_t \comm] \circ \comm \circ [p_1+ t \ha^{p}].
\end{equation}

Soit $\ha^{pa}$ la fonction symétrique avec un paramètre $t$ définie par
\begin{equation}
  \ha^{pa}= p_1 \left( \comm \circ [p_1+t \ha^{p}] \right)
  [1+t \Sigma_t \comm] \circ \comm \circ [p_1+ t \ha^{p}].
\end{equation}

Soit $\ha^{a}$ la fonction symétrique avec un paramètre $t$ définie par
\begin{equation}
  \ha^{a}=[\comm - p_1] \circ [p_1+ t \ha^{p}].
\end{equation}

Soit enfin $\ha$ la fonction symétrique avec un paramètre $t$ définie
par
\begin{equation}
  \ha=\ha^{p}+\ha^{a}-\ha^{pa}.
\end{equation}

Ces relations permettent de calculer ces fonctions symétriques par
récurrence.

On a la relation
\begin{equation}
  \ha^{p}=p_1 \partial_{p_1} \ha.
\end{equation}

\subsection{Caractère des hyperarbres cycliques}

Soit $\hac^{pa}$ la fonction symétrique avec un paramètre $t$ définie par
\begin{equation}
  \hac^{pa}=p_1 [\Sigma_t \assoc] \circ \comm \circ [p_1+ t \hac^{pa}].
\end{equation}

Soit $\hac^{p}$ la fonction symétrique avec un paramètre $t$ définie par
\begin{equation}
  \hac^{p}=p_1 [\Sigma_t \cycle] \circ \comm \circ [p_1+t \hac^{pa}],
\end{equation}
où $\cycle$ est la fonction symétrique dont le terme en degré $n$
correspond à l'action du groupe symétrique sur l'ensemble des ordres
cycliques sur $n$ symboles.

Soit $\hac^{a}$ la fonction symétrique avec un paramètre $t$ définie par
\begin{equation}
  \hac^{a}=[\comm - p_1] \circ [p_1+ t \hac^{pa}].
\end{equation}

Soit enfin $\hac$ la fonction symétrique avec un paramètre $t$ définie
par
\begin{equation}
  \hac=\hac^{p}+\hac^{a}-\hac^{pa}.
\end{equation}

Ces relations permettent de calculer ces fonctions symétriques par
récurrence.

On a la relation
\begin{equation}
  \hac^{p}=p_1 \partial_{p_1} \hac.
\end{equation}

\bibliographystyle{plain}
\bibliography{hypera}

\end{document}